\documentclass[10pt,a4paper]{article}

\usepackage{epsf,epsfig,amsfonts,amsgen,amsmath,amstext,amsbsy,amsopn,amsthm,amssymb
}
\usepackage{ebezier,eepic}
\usepackage{color}
\usepackage{multirow}
\newtheorem{theorem}{Theorem}[section]

\newtheorem{proposition}{Proposition}[section]
\newtheorem{lemma}{Lemma}
\newtheorem{false statement}{False statement}

\theoremstyle{definition}

\newtheorem{claim}{Claim}

\newtheorem{corollary}[claim]{Corollary}
\newtheorem{problem}{Problem}

\begin{document}

\title{\bf\Large Spectral radius and Hamiltonian properties of graphs, II}
\date{}

\author{Jun Ge\thanks{School of Mathematical Sciences,
Sichuan Normal University, Chengdu, 610066, Sichuan, P. R. China.
Email: mathsgejun@163.com}~ and Bo Ning\thanks{Corresponding author.
Center for Applied Mathematics,
Tianjin University, Tianjin, 300072, P. R. China. Email:
bo.ning@tju.edu.cn.}}

\maketitle

\begin{abstract}
In this paper, we first present spectral
conditions for the existence of $C_{n-1}$ in graphs (2-connected
graphs) of order $n$, which are motivated by a conjecture
of Erd\H{o}s. Then we prove spectral conditions for the existence
of Hamilton cycles in balanced bipartite graphs. This result presents
a spectral analog of Moon-Moser's theorem on Hamilton cycles in
balanced bipartite graphs, and extends a previous theorem due to
Li and the second author for $n$ sufficiently large. We conclude
this paper with two problems on tight spectral conditions for the
existence of long cycles of given lengths.
\medskip

\noindent {\bf Keywords:} spectral radius; Hamiltonicity; minimum degree;
long cycle; balanced bipartite graph

\smallskip
\noindent {\bf Mathematics Subject Classification (2010): 05C50, 15A18, 05C38}
\end{abstract}

\section{Introduction}
\begin{figure}[htbp]
\centering
\includegraphics[width=10cm]{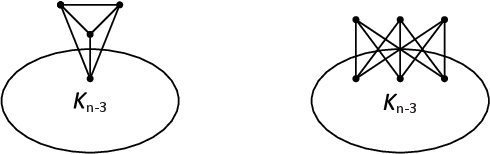}
\renewcommand{\figurename}{Fig.}
\caption{{\footnotesize $L^3_n$ (left) and $N^3_n$ (right).
}}\label{LN}
\end{figure}

Throughout this paper, we only consider graphs which are simple,
finite and undirected. Let $G=(V, E)$ be a graph of order $n$ and
size $e(G)$. Let $S\subset V(G)$. We use $G-S$ to denote the subgraph
induced by $V(G)\backslash S$. If $S$ consists of only one element,
say $S=\{u\}$, then we use $G-u$ instead of $G-\{u\}$.
Let $\lambda_1(G)\geq \lambda_2(G)\geq \cdots \geq \lambda_n(G)$
be all the eigenvalues of the adjacency matrix $A(G)$ of $G$.
Denote by $\lambda(G):=\lambda_1(G)$ the \emph{spectral radius}
of $G$, $q(G)$ the signless Laplacian spectral radius of $G$,
and $\delta(G)$ the \emph{minimum degree} of $G$. Let $G_1$
and $G_2$ be two graphs. We use $G_1+G_2$ to denote the
\emph{disjoint union} of $G_1$ and $G_2$, and $G_1\vee G_2$ to denote
the \emph{join} of $G_1$ and $G_2$. Following some notations
in \cite{LN16}, for $1\leq k\leq(n-1)/2$, we define
$L^k_n=K_1\vee(K_k+K_{n-k-1})$ and $N^k_n=K_k\vee(K_{n-2k}+kK_1)$.
(Fig \ref{LN} illustrates $L^3_n$ and $N^3_n$).
Note that $L^1_{n}=N^1_{n}$. A graph $G$ is called \emph{Hamiltonian}
if it contains a spanning cycle, and is called \emph{pancyclic}
if it contains cycles of lengths from 3 to $v(G)$.
The \emph{circumference} of a graph refers to the length of a
longest cycle in the graph. For terminology
and notations not defined here, we refer the reader to Bondy and
Murty \cite{BM08}.

Many graph theorists have
investigated the relationship between the existence of Hamilton
cycles and paths in graphs and the eigenvalues of some associated
matrices of graphs, for example, see
\cite{BC10,FN10,Z10,NG15,LN16,N16,NL16,LN17,LNP}.
Among these results, the following one has received much attention, which is a corollary of a theorem of Ore \cite{O61} and
of Bondy \cite{B72}, independently.

\begin{theorem}[Fiedler and Nikiforov \cite{FN10}]\label{ThFN}
Let $G$ be a graph of order $n$. If $\lambda(G)>n-2$, then $G$
is Hamiltonian unless $G=N^1_{n}$.
\end{theorem}

By introducing the minimum degree of a graph as a new parameter,
Li and Ning \cite{LN16} extended Theorem \ref{ThFN} in some sense and obtained spectral
analogs of a classical theorem of Erd\H{o}s \cite{E62}.
The following theorem is one of the central results in \cite{LN16}:
Let $k\geq 1$ and $G$ be a graph of order $n\geq\max\{6k+5,(k^2+6k+4)/2\}$.
If $\delta(G)\geq k$ and $\lambda(G)\geq \lambda(N_n^k)$,
then $G$ is Hamiltonian unless $G=N_n^k$.

Since $K_{n-k}\subset N^{k}_n$
and $K_{n-k}\subset L^{k}_n$, we have
$\lambda(L^{k}_n)>n-k-1$ and $\lambda(N^k_n)>n-k-1$.
Nikiforov \cite{N16} further strengthened Li and
Ning's theorem for a graph $G$ of order $n\geq k^3+O(k)$ and
$k\geq 2$, by providing a weaker condition that $\lambda(G)\geq n-k-1$.
Later, the authors \cite{GN1} sharpened the result
mentioned above for the order of graphs by almost a
half.\footnote{For comments on this fact, see Chen and Zhang \cite{CZ18}.}

Since the spectral conditions for $C_n$ are extensively studied, one
can naturally consider similar problems for the possible second longest cycle,
that is, $C_{n-1}$.  Indeed, our first part is motivated by a conjecture
of Erd\H{o}s that says every graph of order $n$ has
a $C_{n-1}$ if its size is at least $\binom{n-2}{2}+4$, which was
confirmed by Bondy \cite{B70}. The edge number condition above is tight, since one can see
the graph $K_1\vee (K_2+K_{n-3})$ has $\binom{n-2}{2}+3$ edges but contains no $C_{n-1}$.

\begin{figure}[htbp]
\centering
\includegraphics[width=4.5cm]{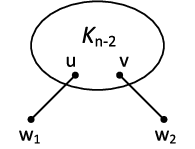}
\renewcommand{\figurename}{Fig.}
\caption{{\footnotesize The graph $\Lambda$.
}}\label{Lambda}
\end{figure}

\begin{theorem}\label{Th(n-1)cycleAQ}
Let $G$ be a graph of order $n\geq 15$.\\
(1) If $\lambda(G)> n-3$,
then $C_{n-1}\subseteq G$, unless $G\subseteq K_1\vee (K_{n-3}+K_2)$ or $G\subseteq \Lambda$
(see Fig \ref{Lambda}).\\
(2) If $q(G)> 2n-6$,
then $C_{n-1}\subseteq G$, unless $G\subseteq K_1\vee (K_{n-3}+K_2)$ or $G\subseteq \Lambda$.
\end{theorem}

By Lemma \ref{compare} in Section \ref{spectral},
Theorem \ref{Th(n-1)cycleAQ} implies the following corollary immediately.

\begin{corollary}\label{Cor1}
Let $G$ be a graph of order $n\geq 15$.\\
(1) If $\lambda(G)\geq\lambda(K_1\vee (K_{n-3}+K_2))$,
then $C_{n-1}\subseteq G$ unless $G=K_1\vee (K_{n-3}+K_2)$.\\
(2) If $q(G)\geq q(K_1\vee (K_{n-3}+K_2))$,
then $C_{n-1}\subseteq G$ unless $G=K_1\vee (K_{n-3}+K_2)$.
\end{corollary}

Considering that the extremal graphs in Corollary 1 contain a cut-vertex,
we then consider similar problems among 2-connected graphs. The answers are as follows.

\begin{theorem}\label{Th2-connected(n-1)cycleAQ}
Let $G$ be a 2-connected graph of order $n\geq 22$.\\
(1) If $\lambda(G)>n-4$,
then $C_{n-1}\subseteq G$ unless $G\subseteq K_2\vee (K_{n-5}+3K_1)$.\\
(2) If $q(G)>2n-8$,
then $C_{n-1}\subseteq G$ unless $G\subseteq K_2\vee (K_{n-5}+3K_1)$.
\end{theorem}

The following corollary follows immediately.

\begin{corollary}
Let $G$ be a 2-connected graph of order $n\geq 22$.\\
(1) If $\lambda(G)\geq\lambda(K_2\vee (K_{n-5}+3K_1))$,
then $C_{n-1}\subseteq G$ unless $G=K_2\vee (K_{n-5}+3K_1)$.\\
(2) If $q(G)\geq q(K_2\vee (K_{n-5}+3K_1))$,
then $C_{n-1}\subseteq G$ unless $G=K_2\vee (K_{n-5}+3K_1)$.
\end{corollary}

In this paper, we also consider Hamilton cycles in
balanced bipartite graphs. Here, a bipartite graph is called
\emph{balanced} if its two partite sets $X$ and $Y$ have
the same number of vertices. Denote by $B^k_n$ the graph
obtained from $K_{n,n}$ by deleting a $K_{k,n-k}$, where $n\geq 2k+1$.
We will improve another theorem of Li and Ning \cite{LN16} on
the spectral condition for Hamilton cycles in balanced
bipartite graphs.
\begin{theorem}\label{balbipHC}
Let $k\geq 1$. Let $G$ be a balanced bipartite graph of order $2n$ with
$\delta(G)\geq k$, where $n\geq k^3+2k+4$.\\
(1) If $G$ is a proper subgraph of $B^k_n$, then $\lambda(G)<\sqrt{n(n-k)}$.\\
(2) If $\lambda(G)\geq \sqrt{n(n-k)}$, then $G$ is Hamiltonian
unless $G=B^k_n$.
\end{theorem}

Since $\lambda(B^k_n)>\sqrt{n(n-k)}$, Theorem \ref{balbipHC} extends the following
theorem in \cite{LN16} due to Li and the second author
(for $n$ sufficiently large).

\begin{theorem}[Li and Ning \cite{LN16}]
Let $k\geq 1$. Let $G$ be a balanced bipartite graph of order $2n$, where $n\geq (k+1)^2$.
If $\delta(G)\geq k$ and $\lambda(G)\geq\lambda(B^k_n)$, then $G$ is
Hamiltonian unless $G=B^k_n$.
\end{theorem}

We organize this paper as follows. In Section 2, we list necessary
preliminaries and prove two structural lemmas.
In Section 3, we introduce the Kelmans operation and list
several spectral inequalities which will be used in
the proof of main theorems. In Section 4,
we prove our main results. We conclude this
paper in the final section with two problems and
some discussions.

\section{Structural results}

In this section, we list several structural theorems that we rely on.

\begin{theorem}[Bondy]\label{ThBondy}
(i) {\rm \cite{B71-J}} Let $G$ be a Hamiltonian graph of order $n$. If $e(G)\geq \frac{n^2}{4}+1$,
then $G$ is pancyclic.
(ii) {\rm \cite{B71-D}} Every graph of order $n\geq 4$ and size $e(G)\geq \binom{n-2}{2}+4$
contains a $C_{n-1}$.
\end{theorem}

The following theorem is a corollary of a theorem of Erd\H{o}s \cite{E62}.
\begin{theorem}[Erd\H{o}s \cite{E62}]\label{Thm_E}
Let $G$ be a 2-connected graph of order $n\geq 13$. If $e(G)\geq \binom{n-2}{2}+5$,
then $G$ is Hamiltonian.
\end{theorem}

The next two theorems on Hamiltonian properties of graphs are also useful for us.
\begin{theorem}[Ore \cite{O63}]\label{Thm_O}
Let $G$ be a graph of order $n$. If $e(G)\geq \binom{n-1}{2}+3$, then $G$ is Hamiltonian-connected.
\end{theorem}

\begin{theorem}[Li and Ning \cite{LN16}]\label{Thm_LN}
Let $G$ be a graph of order $n\geq 6k+5$, where $k\geq 1$. If
$\delta(G)\geq k$ and
$$e(G)>\binom{n-k-1}{2}+(k+1)^2,$$ then $G$ is Hamiltonian unless
$G\subseteq L^k_n$ or $G\subseteq N^k_n$.
\end{theorem}

The following result was proved by Kopylov in \cite{K77} (see the last part of \cite{K77}),
which was originally conjectured by Woodall in \cite{W76}.
\begin{theorem}[Kopylov \cite{K77}]\label{Thkopylov}
Let $n\geq c\geq 5$ and $G$ be a 2-connected graph of order $n$
with circumference less than $c$. If the minimum degree $\delta(G)\geq k\geq 2$,
then
$$e(G)\leq\max\left\{f(n,k,c),f\Big(n,\Big\lfloor\frac{c-1}{2}\Big\rfloor,c\Big)\right\},$$
where $f(n,k,c)=\binom{c-k}{2}+k(n-c+k)$.
\end{theorem}

The following two lemmas, which refine Theorem \ref{ThBondy}(ii),
will play the central role in proving Theorems \ref{Th(n-1)cycleAQ}
and \ref{Th2-connected(n-1)cycleAQ}.
\begin{lemma}\label{Lem1}
Let $G$ be a graph of order $n\geq 15$ and size $e(G)\geq \binom{n-2}{2}$.
Then $G$ contains a $C_{n-1}$, unless $G\subseteq K_1\vee (K_{n-3}+K_2)$
or $G\subseteq \Lambda$.
\end{lemma}
\begin{proof}
We prove the lemma by contradiction. Suppose that $G$ contains no $C_{n-1}$ and
$G\nsubseteq K_1\vee (K_{n-3}+K_2)$,
and $G\nsubseteq \Lambda$.

Suppose that $G$ is Hamiltonian. Since
$e(G)\geq \frac{(n-2)(n-3)}{2}\geq \frac{n^2}{4}+1$ when $n\geq 10$, by Theorem \ref{ThBondy}(i),
$G$ is pancyclic, and thus $C_{n-1}\subseteq G$, a contradiction. Hence $G$ contains no $C_n$ or $C_{n-1}$.
So the circumference of $G$ is less than $n-1$.

Suppose that $G$ is 2-connected. By Theorem \ref{Thkopylov}, $e(G)\leq \max\{f(n,2,n-1),f(n,\lfloor\frac{n}{2}\rfloor-1,n-1)\}$.
Since $f(n,2,n-1)=\binom{n-3}{2}+6$ and $f(n,\lfloor\frac{n}{2}\rfloor-1,n-1)=
\binom{\lceil\frac{n}{2}\rceil}{2}+(\lfloor\frac{n}{2}\rfloor-1)\lfloor\frac{n}{2}\rfloor$,
it is easy to check by WolframAlpha (http://www.wolframalpha.com/) that
$$\max\{f(n,2,n-1),f(n,\lfloor\frac{n}{2}\rfloor-1,n-1)\}=f(n,2,n-1)=\binom{n-3}{2}+6$$ when $n\geq 15$.
However, we have $e(G)\geq \binom{n-2}{2}>\binom{n-3}{2}+6$ when $n\geq 10$, a contradiction. So
$G$ is not 2-connected.

Suppose that $G$ is disconnected. Let $G=G_1\cup G_2$, where $G_1\cap G_2=\emptyset$ and
$v(G_1)\geq v(G_2)\geq 1$. Set $v(G_1)=a\geq \frac{n}{2}$ and $v(G_2)=n-a$.
Obviously, $e(G)\leq \binom{a}{2}+\binom{n-a}{2}$,
which implies that $\binom{a}{2}+\binom{n-a}{2}\geq \binom{n-2}{2}$.
That is, $(a-2)(a+2-n)\geq -1$.
Since $a\geq \frac{n}{2}>7$, we get $a\geq n-2$. Hence $(a,b)=(n-2,2)$ or $(n-1,1)$.

If $(a,b)=(n-2,2)$, then $G\subseteq K_{n-2}+K_2\subset K_1\vee (K_{n-3}+K_2)$, a contradiction.

If $(a,b)=(n-1,1)$, then let us consider $G_1$. Notice that,
$e(G_1)=e(G)\geq \binom{n-2}{2}=\binom{v(G_1)-1}{2}> \binom{v(G_1)-2}{2}+4$
when $v(G_1)\geq 11$ ($n\geq 12$).
By Theorem \ref{Thm_LN},
$G_1$ is either Hamiltonian, which contradicts the fact that $G$ contains no $C_{n-1}$;
or $G_1\subseteq K_{1}\vee (K_{n-3}+K_1)$, which follows
$G\subseteq (K_{1}\vee (K_{n-3}+K_1))+K_1\subset K_1\vee (K_{n-3}+K_2)$,
also a contradiction.

Finally, consider the case that $G$ is connected with a cut-vertex, say $v$.
Let $G=G_1\cup G_2$, where $V(G_1)\cap V(G_2)=\{v\}$ and $v(G_1)\geq v(G_2)\geq 2$.
Set $v(G_1)=a$ and $v(G_2)=b=n+1-a$, where $\frac{n+1}{2}\leq a\leq n-1$. Thus, we have
$\binom{a}{2}+\binom{n+1-a}{2}\geq \binom{n-2}{2}$,
which implies $a^2-(n+1)a+3n-3\geq 0$. That is,
$(a-3)(a-n+2)\geq -3$. If $a<n-2$, then
$(a-3)(a-n+2)\leq 3-a\leq 3-\frac{n+1}{2}< -3$ when $n\geq 12$,
a contradiction. Thus, $n-2\leq a\leq n-1$, which implies $(a,b)=(n-2,3)$ or $(n-1,2)$.

If $(a,b)=(n-2,3)$, then $G\subseteq K_1\vee (K_{n-3}+K_2)$, a contradiction.

If $(a,b)=(n-1,2)$, then $e(G_1)=e(G)-1\geq \binom{n-2}{2}-1=\binom{v(G_1)-1}{2}-1$.
Suppose that $G_1$ is 2-connected. Since
$e(G_1)\geq \binom{v(G_1)-1}{2}-1\geq \binom{v(G_1)-2}{2}+5$ for $v(G_1)\geq 13$ ($n\geq 14$),
by Theorem \ref{Thm_E}, $G_1$ is Hamiltonian. Thus, $C_{n-1}\subseteq G_1\subseteq G$,
a contradiction. So $G_1$ contains a cut-vertex.
Let $G_1=G_{11}\cup G_{12}$, where $V(G_{11})\cap V(G_{12})=\{w\}$ and $v(G_{11})\geq v(G_{12})\geq 2$.
Set $v(G_{11})=s$ and $v(G_{12})=n-s$, where $\frac{n}{2}\leq s\leq n-2$. Thus, we have
$\binom{s}{2}+\binom{n-s}{2}\geq \binom{n-2}{2}-1$,
which implies $s\geq n-2-\frac{4}{s-2}$. Since $s\geq \frac{n}{2}> 7$,
we have $s\geq n-2-\frac{4}{5}$. It follows that $s\geq n-2$.
So $v(G_{11})=n-2$ and $v(G_{12})=2$. Set $V(G_2)=\{v,x\}$ and $V(G_{12})=\{w,y\}$.
If $v=w$, then $G\subseteq K_1\vee (K_{n-3}+2K_1)\subset K_1\vee (K_{n-3}+K_2)$,
a contradiction. Thus, $v\neq w$. If $v=y$, then $G\subseteq K_1\vee (K_{n-3}+2K_1)$,
a contradiction. So $\{x,v\}\cap \{y,w\}=\emptyset$. In this case,
$G\subseteq\Lambda$.

The proof of Lemma \ref{Lem1} is complete.
\end{proof}

\begin{lemma}\label{Lem2}
Let $G$ be a 2-connected graph of order $n\geq 22$ and size $e(G)\geq \binom{n-3}{2}-2$.
Then $G$ contains a $C_{n-1}$, unless $G\subseteq K_2\vee (K_{n-5}+3K_1)$.
\end{lemma}

\begin{proof}
We prove the lemma by contradiction. Suppose $G$ contains no $C_{n-1}$,
and $G\nsubseteq K_2\vee (K_{n-5}+3K_1)$.
We shall first prove two claims.
\setcounter{claim}{0}
\begin{claim}\label{claim-notHamiltonian}
$G$ is not Hamiltonian.
\end{claim}
\begin{proof}
Suppose that $G$ is Hamiltonian. Notice that we have
$e(G)\geq \binom{n-3}{2}-2\geq\frac{n^2}{4}+1$ when $n\geq 14$. By
Theorem \ref{ThBondy}(i), $G$ is pancyclic, and thus contains
a $C_{n-1}$, a contradiction.
\end{proof}

\begin{claim}
$G$ contains a 2-cut.
\end{claim}
\begin{proof}
Suppose that $G$ is 3-connected. Then $\delta(G)\geq 3$. By Theorem \ref{Thkopylov},
we have
$$e(G)\leq \max\{f(n,3,n-1),f(n, \lfloor\frac{n}{2}\rfloor-1, n-1)\},$$
where $f(n,k,c)=\binom{c-k}{2}+k(n-c+k)$.

Since $f(n,3,n-1)=\binom{n-4}{2}+12$ and
$f(n,\lfloor\frac{n}{2}\rfloor-1,n-1)=
\binom{\lceil\frac{n}{2}\rceil}{2}+(\lfloor\frac{n}{2}\rfloor-1)\lfloor\frac{n}{2}\rfloor$,
it is easy to check by WolframAlpha (http://www.wolframalpha.com/) that
$$\max\{f(n,3,n-1),f(n,\lfloor\frac{n}{2}\rfloor-1,n-1)\}=f(n,3,n-1)=\binom{n-4}{2}+12$$
when $n\geq 22$.
However, we get $e(G)\geq \binom{n-3}{2}-2>\binom{n-4}{2}+12$ when $n\geq 19$, a contradiction.

Now we know $G$ is 2-connected but not 3-connected, so $G$ contains a 2-cut. This proves the claim.
\end{proof}

We choose $G_1$ and $G_2$ such that:\\
(i) $G=G_1\cup G_2$ with $V(G_1)\cap V(G_2)=\{u,v\}$, where $\{u,v\}$ is a 2-cut; \\
(ii) $v(G_1)-v(G_2)$ is as large as possible.

In the following, we call $\{G_1,G_2\}$ a \emph{good pair}, if it satisfies both (i) and (ii).
Set $v(G_1)=a$ and $v(G_2)=b$. Then $\frac{n+2}{2}\leq a\leq n-1$ and $v(G_2)=n+2-a$. Hence
$\binom{a}{2}+\binom{n+2-a}{2}-1\geq e(G)\geq \binom{n-3}{2}-2$.
Therefore, $a^2-(n+2)a+5n-4\geq 0$, and hence,
\begin{align}\label{Al_2con-1}
(a-5)(a-n+3)\geq -11.
\end{align}
Since $n\geq 22$, $n-4\geq \frac{n+2}{2}$.
If $a\leq n-4$, then
$(a-5)(a-(n-3))\leq ((n-4)-5)((n-4)-(n-3))=-(n-9)<-11$
when $n\geq 21$, a contradiction to (\ref{Al_2con-1}).
Thus $a\geq n-3$, which implies that $(a,b)=(n-3,5)$, or $(a,b)=(n-2,4)$, or $(a,b)=(n-1,3)$.

Note that in the case $(a,b)=(n-3,5)$ or $(a,b)=(n-2,4)$,
there are no vertices of degree two in $G$, otherwise
$\{G_1,G_2\}$ is not a good pair.

Suppose that $(a,b)=(n-3,5)$. Let $V(G_2)\backslash \{u,v\}=\{w_1,w_2,w_3\}$.
We obtain
$$e(G_1)\geq e(G)-e(G_2)\geq \binom{n-3}{2}-2-\binom{5}{2}\geq \binom{n-4}{2}+3=\binom{v(G_1)-1}{2}+3$$
when $n\geq 19$.
By Theorem \ref{Thm_O}, $G_1$ is Hamiltonian-connected.
Take any $(u,v)$-Hamilton path in $G_1$, say $P_1$. We can see
$G_2-\{u,v\}$ is connected, otherwise there is a vertex of degree 2.
If $G_2-\{u,v\}\cong K_3$, since $G$ is 2-connected,
there exists $(u,v)$-Hamilton path in $G_2$,
say $P_2$. Then $P_1\cup P_2$ is a Hamilton cycle in $G$, a contradiction to
Claim \ref{claim-notHamiltonian}. If $G_2-\{u,v\}\cong P_3$,
then suppose $G_2-\{u,v\}$ is the path $w_1w_2w_3$.
Since there is no vertex of degree 2 in $G$,
$u$ and $v$ are neighbours of both $w_1$ and $w_3$.
Then $vP_1uw_1w_2w_3v$ is a Hamilton cycle in $G$, a contradiction to
Claim \ref{claim-notHamiltonian}.

Suppose that $(a,b)=(n-2,4)$. Let $V(G_2)\backslash \{u,v\}=\{w_1,w_2\}$.
Since there is no vertex of degree 2 in $G$,
$\{w_1u, w_1v, w_2u, w_2v, w_1w_2\}\subset E(G)$.
Let $H:=G-w_1$.
Then $\delta(H)\geq 2$. We obtain
$$e(H)=e(G)-3\geq \binom{n-3}{2}-5>\binom{n-4}{2}+9=\binom{v(H)-3}{2}+9$$
for $n\geq 19$.
By Theorem \ref{Thm_LN}, $H$ is Hamiltonian, unless $H\subseteq L^2_{n-1}$ or
$H\subseteq N^2_{n-1}$. If $H$ is Hamiltonian,
then $C_{n-1}\subset H\subset G$, a contradiction.
Next, we shall show that
$H\nsubseteq L^2_{n-1}$.
Since $H=G-w_1$ and $N_G(w_1)=\{u, v, w_2\}$,
it is easy to prove $H$ is 2-connected.
But $L^2_{n-1}$ contains a cut-vertex.
Hence $H\nsubseteq L^2_{n-1}$, and it follows
$H\subseteq N^2_{n-1}$. Let $\{x_1,x_2\}\subset V(H)$ such that $x_1,x_2$ are two
vertices of degree 2 in $N^2_{n-1}$.
Since there is no vertex of degree 2 in $G$,
$x_1$ and $x_2$ must be neighbours of $w_1$.
Since $N_G(w_1)=\{u, v, w_2\}$, $\{x_1, x_2\}\cap \{u, v\}\neq\emptyset$.
Hence at least one of $u, v$ has degree 3, say $u$.
Furthermore, the neighbour of $u$ in $G$ other than $w_1$ and $w_2$ must not be $v$,
otherwise $v$ is a cut-vertex in $G$, which contradicts the fact $G$ is 2-connected.
Let the neighbour of $u$ in $G$ other than $w_1$ and $w_2$ be $z$.
Since
$$e(G-\{u,w_1,w_2\})=e(G)-6\geq \binom{n-3}{2}-8\geq \binom{n-4}{2}+3$$
for $n\geq 15$, we obtain that $G-\{u,w_1,w_2\}$ is Hamiltonian-connected
by Theorem \ref{Thm_O}. Take a $(z,v)$-Hamilton path $P$ in
$G-\{u,w_1,w_2\}$. We can see $vPzuw_1w_2v$ is a Hamilton cycle in $G$,
a contradiction to Claim \ref{claim-notHamiltonian}.

\begin{figure}[htbp]
\centering
\includegraphics[width=13cm]{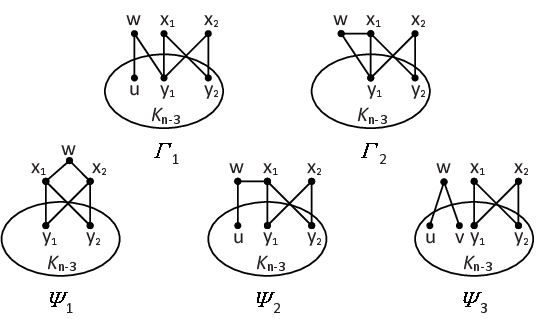}
\renewcommand{\figurename}{Fig.}
\caption{{\footnotesize $\Gamma_1$, $\Gamma_2$, $\Psi_1$, $\Psi_2$ and $\Psi_3$.
}}\label{gamma}
\end{figure}

Finally, consider the case that $(a,b)=(n-1,3)$. Let $V(G_2)=\{u,v,w\}$
and $H:=(G-w)\cup \{uv\}$. Then $H$ is 2-connected.
Moreover, $e(H)\geq e(G)-d_G(w)\geq \binom{n-3}{2}-4$. When $n\geq 18$, we have
$e(H)>\binom{v(H)-3}{2}+9=\binom{n-4}{2}+9$.
By Theorem \ref{Thm_LN}, when $v(H)\geq 17$ (that is, $n\geq 18)$,
$H$ is Hamiltonian unless $H\subseteq L^2_{n-1}$ or $H\subseteq N^2_{n-1}$.
Since $H$ is 2-connected and $L^2_{n-1}$ has a cut-vertex, $H\nsubseteq L^2_{n-2}$.
Hence $H$ is Hamiltonian or $H\subseteq N^2_{n-1}$.
Assume that $H$ is Hamiltonian.
If $uv\in E(G)$, then $G$ contains a $C_{n-1}$, a contradiction. Thus, $uv\notin E(G)$.
If the Hamilton cycle in $H$, say $C$, does not pass through the edge $uv$, then it
is also in $G$, a contradiction. Thus, $C$ passes through $uv$, and hence
there is a $(u,v)$-Hamilton path in $G_1$. Together with the path $uwv$,
we can find a Hamilton cycle in $G$, a contradiction. So
$H\subseteq N^2_{n-1}=K_2\vee(K_{n-5}+2K_1)$.
Let $\{x_1,x_2\}\subseteq V(N^2_{n-1})$ such that each is of degree 2 and
let $\{y_1,y_2\}$ be the 2-cut in $N^2_{n-1}$.
Set $V(H)=V(N^2_{n-1})$. Since $H$ is 2-connected, $y_1,y_2$ are
still neighbors of $x_1,x_2$ in $H$.
So, according to the locations of $y_1,y_2$, we obtain the following subcases
(in the sense of isomorphism):
(1) If $\{u,v\}=\{y_1,y_2\}$, then $G\subseteq K_2\vee (K_{n-5}+3K_1)$, a contradiction;
(2) If $|\{u,v\}\cap \{y_1,y_2\}|=1$, then $G$ is a subgraph of $\Gamma_1$ or $\Gamma_2$
(see Fig \ref{gamma});
(3) If $\{u,v\}\cap \{y_1,y_2\}=\emptyset$,
then $G$ is a subgraph of a graph in $\Psi_1$, $\Psi_2$, or $\Psi_3$
(see Fig \ref{gamma}). By Proposition \ref{exceptions}
(whose proof will be presented later), both subcases (2) and (3) contain either
a $C_{n-1}$ or a $C_n$, also a contradiction.

The proof is complete.
\end{proof}

\begin{proposition}\label{exceptions}
Let $G$ be a graph of order $n\geq 21$ and size $e(G)\geq \binom{n-3}{2}-2$.
Suppose that the degrees of $w,x_1,x_2$ in $G$ equal
the degrees of those in $\Gamma_i$ for $i=1,2$
and in $\Psi_j$ for $j=1,2,3$, respectively.
If $G$ is a spanning subgraph of $\Psi_1$ or $\Psi_2$, then $G$ is Hamiltonian.
If $G$ is a spanning subgraph of $\Gamma_1$, $\Gamma_2$ or $\Psi_3$,
then $G$ contains a $C_{n-1}$.
\end{proposition}

\begin{proof}
(1) If $G\subseteq \Psi_1$, then let us consider $G-\{w,x_1,x_2\}$.
Since $v(G-\{w,x_1,x_2\})=n-3$ and
$$e(G-\{w,x_1,x_2\})=e(G)-6\geq \binom{n-3}{2}-8\geq \binom{n-4}{2}+3$$
when $n\geq 15$, we obtain that $G-\{w,x_1,x_2\}$ is Hamiltonian-connected
by Theorem \ref{Thm_O}. So there is a $(y_1, y_2)$-Hamilton path $P$ in $G-\{w,x_1,x_2\}$.
Therefore $y_2Py_1x_1wx_2y_2$ is a Hamilton cycle in $G$.

(2) If $G\subseteq \Psi_2$, then let us consider $G-\{w,x_1,x_2,y_1\}$.
Since $v(G-\{w,x_1,x_2,y_1\})=n-4$ and
$$e(G-\{w,x_1,x_2,y_1\})=e(G)-4-d_G(y_1)\geq \binom{n-3}{2}-2-4-(n-2)\geq \binom{n-5}{2}+3$$
when $n\geq 16$, we obtain that $G-\{w,x_1,x_2,y_1\}$ is Hamiltonian-connected
by Theorem \ref{Thm_O}. So there is a $(u, y_2)$-Hamilton path $P$ in $G-\{w,x_1,x_2,y_1\}$.
Therefore, $y_2Puwx_1y_1x_2y_2$ is a Hamilton cycle in $G$.

(3) If $G\subseteq \Gamma_1$, we consider $G-\{w,x_1,x_2,y_1\}$.
Since $v(G-\{w,x_1,x_2,y_1\})=n-4$ and
$$e(G-\{w,x_1,x_2,y_1\})=e(G)-3-d_G(y_1)\geq \binom{n-3}{2}-2-3-(n-1)\geq \binom{n-5}{2}+3$$
when $n\geq 16$, we obtain that $G-\{w,x_1,x_2,y_1\}$ is Hamiltonian-connected
by Theorem \ref{Thm_O}. So there is a $(u, y_2)$-Hamilton path $P$ in $G-\{w,x_1,x_2,y_1\}$.
Therefore, $y_2Puwy_1x_1y_2$ is a $C_{n-1}$ in $G$.

(4) If $G\subseteq \Gamma_2$, we consider $G-\{w,x_1,x_2\}$.
Since $v(G-\{w,x_1,x_2\})=n-3$ and
$$e(G-\{w,x_1,x_2\})=e(G)-6\geq \binom{n-3}{2}-8\geq \binom{n-4}{2}+3$$
when $n\geq 15$, we obtain that $G-\{w,x_1,x_2\}$ is Hamiltonian-connected
by Theorem \ref{Thm_O}. So there is a $(y_1, y_2)$-Hamilton path $P$ in $G-\{w,x_1,x_2\}$.
Therefore, $y_2Py_1wx_1y_2$ is a $C_{n-1}$ in $G$.

(5) If $G\subseteq \Psi_3$, then we shall prove that there is a $C_{n-1}$ in $G$.
First we claim that $N(v)\cap N(y_1)\neq \emptyset$, otherwise $d_G(v)+d_G(y_1)\leq n$,
and it follows that
$$e(G)\leq e(G-\{w,x_1,x_2,v,y_1\})+d_G(v)+d_G(y_1)+3\leq \binom{n-5}{2}+n+3<\binom{n-3}{2}-2$$
when $n>14$, a contradiction. Hence there exists a vertex $z\in N(v)\cap N(y_1)$.
We consider the graph $G-\{w,x_1,x_2,v,y_1,z\}$.
Since $v(G-\{w,x_1,x_2,v,y_1,z\})=n-6$ and
$$e(G-\{w,x_1,x_2,v,y_1,z\})\geq e(G)-6-3(n-4)\geq \binom{n-3}{2}-3n+4\geq \binom{n-7}{2}+3$$
when $n\geq 21$, we obtain that $G-\{w,x_1,x_2,v,y_1,z\}$ is Hamiltonian-connected
by Theorem \ref{Thm_O}. So there is a $(u, y_2)$-Hamilton path $P$ in $G-\{w,x_1,x_2,v,y_1,z\}$.
Therefore, $y_2Puwvzy_1x_1y_2$ is a $C_{n-1}$ in $G$.

The proof is complete.
\end{proof}

The following theorem will be used in the proof of Theorem \ref{balbipHC}.

\begin{theorem}[Li and Ning \cite{LN16}]\label{ThLN16-BalBip}
Let $G$ be a balanced bipartite graph of order $2n$, where $n\geq 2k+1$, $k\geq 1$. If
$\delta(G)\geq k$ and $e(G)>n(n-k-1)+(k+1)^2$, then $G$ is Hamiltonian unless $G\subseteq B_n^k$.
\end{theorem}

\section{Spectral inequalities}\label{spectral}
First, we introduce the Kelmans operation \cite{K81}.
Given a graph $G$ and two specified vertices $u$ and $v$,
we construct a new graph $G[u\rightarrow v]$ as follows:
we delete all edges between $u$ and $S:= N(u)\backslash (N(v) \cup \{v\})$
and add all edges between $v$ and $S$. In notation, $G[u\rightarrow v]$ is defined as:
$$V(G[u\rightarrow v])=V(G) ~\text{and}~
E(G[u\rightarrow v])=(E(G)\backslash\{uw:w\in S\})\cup \{vw:w\in S\}.$$
We call it \emph{the Kelmans operation} (from $u$ to $v$).

Csikv\'ari \cite{C09} proved that the Kelmans operation does not decrease the spectral radius of a graph.

\begin{theorem}[Csikv\'ari \cite{C09}]\label{Kelmans}
Let $G$ be a graph and $u,v$ be two vertices of $G$.
Let $G'=G[u\rightarrow v]$.
Then $\lambda(G')\geq \lambda(G)$.
\end{theorem}

Li and Ning \cite{LN16} proved a signless spectral radius version.

\begin{theorem}[Li and Ning \cite{LN16}]\label{q-Kelmans}
Let $G$ be a graph and $u,v$ be two vertices of $G$.
Let $G'=G[u\rightarrow v]$.
Then $q(G')\geq q(G)$.
\end{theorem}

Generally speaking, we use these two theorems
to determine the extremal graphs, if we already
have turned the original problems into similar ones under the condition
of number of edges.

The next result helps us to determine the extremal graphs with the help of the two theorems above.
\begin{lemma}\label{compare}
$\lambda(K_1 \vee (K_{n-3}+K_2))>\lambda(\Lambda)$ and $q(K_1 \vee (K_{n-3}+K_2))>q(\Lambda)$.
\end{lemma}

\begin{proof}
Let $G':=\Lambda[u\rightarrow v]$ (see Fig \ref{Lambda} for the
vertices $u$, $v$ and the graph $\Lambda$). Then $G'=K_1 \vee (K_{n-3}+2K_1)$.
By Theorems \ref{Kelmans} and \ref{q-Kelmans}, we obtain
$\lambda(K_1 \vee (K_{n-3}+2K_1))\geq \lambda(\Lambda)$ and
$q(K_1 \vee (K_{n-3}+2K_1))\geq q(\Lambda)$, respectively. Since
$K_1 \vee (K_{n-3}+2K_1)\subset K_1 \vee (K_{n-3}+K_2)$ and $K_1 \vee (K_{n-3}+K_2)$ is connected,
Lemma \ref{compare} is proved.
\end{proof}

Finally, we need several inequalities below on spectral radius
or signless spectral radius in terms of number of edges and vertices
for graphs or bipartite graphs. We mainly use them to get
a sufficient condition in terms of number of edges for each problem.

\begin{theorem}[Nosal \cite{N70}, Bhattacharya, Friedland, and Peled \cite{BFP08}]\label{ThBFP08}
Let $G$ be a bipartite graph. Then
$\lambda(G)\leq \sqrt{e(G)}$.
\end{theorem}

\begin{theorem}[Hong \cite{H93}]\label{ThHong}
Let $G$ be a graph of order $n$. If $\delta(G)\geq 1$, then
$$\lambda(G)\leq\sqrt{2e(G)-n+1}.$$
\end{theorem}

\begin{theorem}[Hong, Shu, and Fang \cite{HSK00}]\label{ThHSF}
Let $G$ be a connected graph of order $n$ and size $e(G)$.
If minimum degree $\delta(G)\geq k\geq 1$, then
$$\lambda(G)\leq \frac{k-1+\sqrt{(k+1)^2+4(2e(G)-kn)}}{2}.$$
\end{theorem}

\begin{theorem}[Feng and Yu \cite{FY09}]\label{ThFY}
Let $G$ be a graph of order $n$. Then
$$q(G)\leq \frac{2e(G)}{n-1}+n-2.$$
\end{theorem}

\section{Proofs}

\noindent
{\bf Proof of Theorem \ref{Th(n-1)cycleAQ}(1).}
Suppose that $\delta(G)\geq 1$. Then by Theorem \ref{ThHong} and the assumption, we obtain
$$\sqrt{2e(G)-n+1}\geq \lambda(G)> n-3,$$
which implies that $2e(G)> (n-2)(n-3)+2$.
That is,
\begin{align}
e(G)\geq\binom{n-2}{2}+2.
\end{align}
By Lemma \ref{Lem1}, $G$ contains a $C_{n-1}$, or $G\subseteq K_1\vee (K_{n-3}+K_2)$,
or $G\subseteq \Lambda$.

Now we turn to the case that $G$ contains isolated vertices. Since
the maximum degree $\triangle(G)\geq \lambda(G)> n-3$,
we have $\triangle(G)\geq n-2$, which follows $G$ contains exactly one isolated vertex, say $u$.
Let $H=G-u$. Then $\delta(H)\geq 1$ and $v(H)=n-1$.
Again, by Theorem \ref{ThHong}, we obtain $\sqrt{2e(H)-v(H)+1}=\sqrt{2e(G)-n+2}> n-3$.
We get $2e(G)> n^2-5n+7=(n-2)(n-3)+1$.
Thus, $e(G)\geq \binom{n-2}{2}+1$.
By Lemma \ref{Lem1}, $G$ contains a $C_{n-1}$,
or $G\subseteq K_1\vee (K_{n-3}+K_2)$, or $G\subseteq \Lambda$.
The proof is complete. {\hfill$\Box$}

\noindent
{\bf Proof of Theorem \ref{Th(n-1)cycleAQ}(2).}
Recall Theorem \ref{ThFY}. We obtain $\frac{2e(G)}{n-1}+n-2\geq q(G)> 2n-6$,
which follows $2e(G)>(n-1)(n-4)$. Since $2e(G)$ and $(n-1)(n-4)$ are both
even, we deduce that $2e(G)\geq (n-1)(n-4)+2$,
that is, $e(G)\geq \binom{n-2}{2}$.
By Lemma \ref{Lem1}, $G$ contains a $C_{n-1}$,
or $G\subseteq K_1\vee (K_{n-3}+K_2)$, or $G\subseteq \Lambda$.
The proof is complete.  {\hfill$\Box$}

\noindent
{\bf Proof of Theorem \ref{Th2-connected(n-1)cycleAQ}(1).}
Since $G$ is 2-connected, we get $\delta(G)\geq 2$. By Theorem \ref{ThHSF},
we have $\lambda(G)\leq \frac{1+\sqrt{9+4(2e(G)-2n)}}{2}$.
Hence $\frac{1+\sqrt{9+4(2e(G)-2n)}}{2}>n-4$, which follows that
$2e(G)>(n-3)(n-4)+6$. That is,
$e(G)> \binom{n-3}{2}+3$.
By Lemma \ref{Lem2}, $G$ contains a $C_{n-1}$, or $G\subseteq K_2\vee (K_{n-5}+3K_1)$.
The proof is complete.    {\hfill$\Box$}

\noindent
{\bf Proof of Theorem \ref{Th2-connected(n-1)cycleAQ}(2).}
By Theorem \ref{ThFY}, we obtain $\frac{2e(G)}{n-1}+n-2\geq q(G)> 2n-8$,
which follows $2e(G)>(n-1)(n-6)$. Since $2e(G)$ and $(n-1)(n-6)$ are both
even, we infer that $2e(G)\geq (n-1)(n-6)+2$,
that is, $e(G)\geq \binom{n-3}{2}-2$.
By Lemma \ref{Lem2}, $G$ contains a $C_{n-1}$, or $G\subseteq K_2\vee (K_{n-5}+3K_1)$.
The proof is complete.  {\hfill$\Box$}

\begin{figure}[htbp]
\centering
\includegraphics[width=10cm]{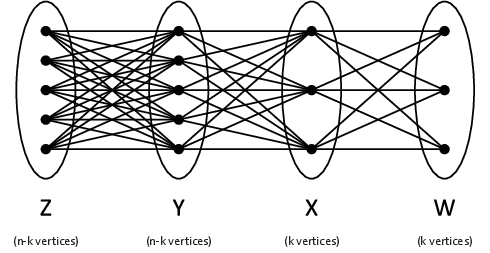}
\renewcommand{\figurename}{Fig.}
\caption{{\footnotesize $X$, $Y$, $Z$, and $W$ in $B_n^k$.
}}\label{B_n^k}
\end{figure}

\noindent
{\bf Proof of Theorem \ref{balbipHC}(1).}  Write $W$ for the set of vertices of $B_n^k$ of degree $k$.
Let $X=N(W)$, $Y=N(X)-W$, and $Z=N(Y)-X$ (see Fig \ref{B_n^k}).
Note that $|W|=|X|=k$ and $|Y|=|Z|=n-k$.

\setcounter{claim}{0}
\begin{claim}
$\lambda(G)\leq \lambda(B_n^k-e)$, where $e=uv$, $u\in Y$, $v\in Z$.
\end{claim}
For any proper subgraph $G$ of $B_n^k$,
since $\delta(G)\geq k$, $G$ must contain all the edges incident to $W$.
Thus,
either $G\subseteq G_1:=B_n^k-e_1$, where $e_1=uv\in E(X,Y)$,
or $G\subseteq G_2:=B_n^k-e_2$, where $e_2=uv\in E(Y,Z)$.
So $\lambda(G)\leq \max\{\lambda(G_1),\lambda(G_2)\}$

Recall $\{u, v\}\subset X\cup Y\cup Z$.
By symmetry, there are two cases: (i) $u\in X$, $v\in Y$;
(ii) $u\in Y$, $v\in Z$.

Let $u'$ be a vertex in $X$, $v'$ be a vertex in $Y$, and $w'$ be a vertex in $Z$. Then
$(B_n^k-u'v')[w'\rightarrow u']=B_n^k-v'w'$.
By Theorem \ref{Kelmans}, $\lambda(B_n^k-v'w')\geq \lambda(B_n^k-u'v')$.
We have proved Claim 1.
\begin{claim}
$\lambda(B_n^k-e)<\sqrt{n(n-k)}$, where $e=uv$, $u\in Y$, $v\in Z$.
\end{claim}
We prove the claim by contradiction. Let $G'=B_n^k-e$, where $e=uv$, $u\in Y$, $v\in Z$.
Set $\lambda:=\lambda(G')$. Let $\mathbf{x}=(x_1, \ldots, x_{2n})$
be a positive unit eigenvector to $\lambda$.
Suppose that
\begin{align}
\lambda(G')\geq \sqrt{n(n-k)}. \label{con3}
\end{align}
Let
\begin{align*}
w  &  :=x_{i},\ i\in W,\\
x  &  :=x_{i},\ i\in X,\\
y  &  :=x_{i},\ i\in Y\backslash \{u\},\\
z  &  :=x_{i},\ i\in Z\backslash\{v\}, \\
s  &  :=x_u, \\
t  &  :=x_v.
\end{align*}

Note that the $2n$ eigenequations of $G'$ are reduced to the
following six types:
\begin{align}
\lambda w & =kx, \label{w1}\\
\lambda x & =kw+(n-k-1)y+s, \label{x1}\\
\lambda y & =kx+(n-k-1)z+t, \label{y1}\\
\lambda z & =(n-k-1)y+s, \label{z1}\\
\lambda s & =kx+(n-k-1)z, \label{s1}\\
\lambda t & =(n-k-1)y. \label{t1}
\end{align}
From (\ref{y1}) and (\ref{s1}) we have
\begin{align*}
\lambda y-\lambda s=\left[kx+t+(n-k-1)z\right]-\left[kx+(n-k-1)z\right],
\end{align*}
that is,
\begin{align}
t=\lambda(y-s). \label{t2}
\end{align}
From (\ref{z1}) and (\ref{t1}) we have
\begin{align*}
\lambda z-\lambda t=(n-k-1)y+s-(n-k-1)y,
\end{align*}
that is,
\begin{align}
s=\lambda(z-t). \label{s2}
\end{align}
By putting (\ref{t2}) into (\ref{s2}), we obtain
\begin{align}
s = \lambda\left[z-\lambda(y-s)\right]= \frac{\lambda^2 y-\lambda z}{\lambda^2-1}. \label{s}
\end{align}
Hence
\begin{align}
t = \lambda(y-s) = \lambda\left(y-\frac{\lambda^2 y-\lambda z}{\lambda^2-1}\right)
  = \frac{\lambda^2 z-\lambda y}{\lambda^2-1}. \label{t}
\end{align}
By using (\ref{s}), equation (\ref{z1}) becomes
\begin{align*}
\lambda z=(n-k-1)y+\frac{\lambda^2y-\lambda z}{\lambda^2-1},
\end{align*}
from which it follows that
\begin{align}
z=\frac{\lambda^2-1}{\lambda^3}\left(n-k+\frac{1}{\lambda^2-1}\right)y. \label{e1}
\end{align}
Since $4<n-k=\lambda(K_{n-k, n-k})<\lambda<\lambda(K_{n, n})=n,$
we obtain
$$z<\frac{\lambda^2-1}{\lambda^3}\left(\lambda+\frac{1}{\lambda^2-1}\right)y
=\frac{\lambda^3-\lambda+1}{\lambda^3}\cdot y<y.$$
Let $f(x)=\frac{x^2-1}{x^3}$. Then $f'(x)=\frac{3-x^2}{x^4}<0$ when $x>\sqrt{3}$.
So $\frac{\lambda^2-1}{\lambda^3}$ decreases when $\lambda>\sqrt{3}$,
which follows
\begin{align*}
z & = \frac{\lambda^2-1}{\lambda^3}(n-k)y+\frac{y}{\lambda^3} \\
  & > \frac{n^2-1}{n^3}(n-k)y+\frac{y}{n^3} \\
  & = \frac{n^3-kn^2-n+k+1}{n^3}\cdot y \\
  & > \frac{n^3-kn^2-n}{n^3}\cdot y \\
  & = \left(1-\frac{k}{n}-\frac{1}{n^2}\right)y.
\end{align*}
Therefore,
\begin{align}
\left(1-\frac{k}{n}-\frac{1}{n^2}\right)y<z<y. \label{in1}
\end{align}

Note that if we remove all edges between $W$ and $X$ and add the edge
$uv$ to $G'$, we obtain the graph $K_{n, n-k}+kK_1$. Let
$\mathbf{x}^{\prime\prime}$ be the restriction of $\mathbf{x}$
to $K_{n, n-k}$. We find that
$$\left\langle A(K_{n, n-k})\mathbf{x}^{\prime\prime},\mathbf{x}^{\prime\prime}\right\rangle =\left\langle A(G')\mathbf{x},\mathbf{x}\right\rangle +2st-2k^2wx=\lambda+2(st-k^2wx).$$
Since $\|\mathbf{x}^{\prime\prime}\|<1$,
$\left\langle A(K_{n, n-k})\mathbf{x}^{\prime\prime},\mathbf{x}^{\prime\prime}\right\rangle <\lambda(K_{n, n-k})=\sqrt{n(n-k)}$,
that is,
\begin{align}
\lambda+2(st-k^2wx) < \sqrt{n(n-k)} \label{in2}.
\end{align}

Recall that $\lambda\geq \sqrt{n(n-k)}$ (see (\ref{con3})). This assumption, together
with (\ref{in2}) yields $st-k^2wx<0$.
Recall that $\lambda w=kx$ (see (\ref{w1})). We rewrite it by
\begin{align*}
\lambda(st-k^2wx)=\lambda st-k^3x^2<0,
\end{align*}
that is,
\begin{align}
\lambda st<k^3x^2. \label{in3}
\end{align}
Noting (\ref{in1}), we have
\begin{align*}
\lambda s&=\lambda\cdot\frac{\lambda^2y-\lambda z}{\lambda^2-1}
> \lambda\cdot\frac{\lambda^2y-\lambda y}{\lambda^2-1}
=\frac{\lambda^2}{\lambda+1}\cdot y>(\lambda-1)y>(n-k-1)y\geq(k^3+k+3)y,
\end{align*}
and
\begin{align*}
t&=\frac{\lambda^2z-\lambda y}{\lambda^2-1}>z-\frac{y}{\lambda-1}
>\left[\left(1-\frac{k}{n}-\frac{1}{n^2}\right)y-\frac{y}{n-k-1}\right]
=\left(1-\frac{k}{n}-\frac{1}{n^2}-\frac{1}{n-k-1}\right)y\\
&>\left(1-\frac{k}{n}-\frac{1}{2n}-\frac{3}{2n}\right)y
=\left(1-\frac{k+2}{n}\right)y\geq \left(1-\frac{k+2}{k^3+2k+4}\right)y.
\end{align*}

Thus, we can estimate
the left side of the inequality in (\ref{in3}) as follows:
\begin{eqnarray*}
\ \lambda st & > & (k^3+k+3)\left(1-\frac{k+2}{k^3+2k+4}\right)y^2 \\
& = &\frac{k^3(k^3+2k+5)+(k+2)(k+3)}{k^3+2k+4}\cdot y^2 \\
& > & k^3y^2.
\end{eqnarray*}
Together with (\ref{in3}), we have
\begin{align}
y^2<x^2. \label{in4}
\end{align}

From (\ref{w1}), (\ref{x1}), (\ref{s}) and (\ref{e1}), we have
\begin{eqnarray*}
\left(\lambda-\frac{k^2}{\lambda}\right)x
& = & (n-k-1)y+s \\
& = & (n-k-1)y+\frac{\lambda^2 y-\lambda z}{\lambda^2-1} \\
& = & (n-k)y+\frac{y-\lambda\cdot\frac{\lambda^2-1}{\lambda^3}
\left(n-k+\frac{1}{\lambda^2-1}\right)y}{\lambda^2-1} \\
& = & (n-k)y-\frac{n-k-1}{\lambda^2}y\\
& < &(n-k)y.
\end{eqnarray*}
Since
\begin{align*}
n\geq k^3+2k+4=(k^3+k)+k+4\geq 2\sqrt{k^3\cdot k}+k+4=2k^2+k+4>\frac{3}{2}k^2+\frac{1}{2}k+\frac{1}{24},
\end{align*}
we have
\begin{align*}
n-\frac{k}{2}-\frac{1}{12}
=\sqrt{n^2-kn-\frac{1}{6}\left[n-\left(\frac{3}{2}k^2+\frac{1}{2}k+\frac{1}{24}\right)\right]}
<\sqrt{n(n-k)}.
\end{align*}
Therefore, we obtain
\begin{align}
n-k<n-\frac{k}{2}-\frac{1}{12}<\sqrt{n(n-k)}\leq \lambda <\lambda (K_{n, n})=n. \label{lambda}
\end{align}
Note that $k^3+k+4>3k^2$ for all $k\geq 1$, hence $\lambda>n-k\geq k^3+k+4>3k^2$,
and it follows $\frac{k^2}{\lambda}<\frac{1}{3}$.
Therefore,
$$x^2<\left(\frac{n-k}{\lambda-\frac{k^2}{\lambda}}\right)^2y^2< \left(\frac{n-k}{n-\frac{k}{2}-\frac{1}{12}-\frac{1}{3}}\right)^2y^2< y^2,$$
contradicting (\ref{in4}).
Now we have proved Claim 2.

Together with Claims 1 and 2, the proof is complete.{\hfill$\Box$}

\noindent
{\bf Proof of Theorem \ref{balbipHC}(2)}
By the initial condition and Theorem \ref{ThBFP08},
$\sqrt{n(n-k)}\leq\lambda(G)\leq\sqrt{e(G)}$.
Thus, we obtain $$e(G)\geq n(n-k)>n(n-k-1)+(k+1)^2$$
when $n\geq (k+1)^2+1$. Notice that $k^3+2k+2\geq (k+1)^2+1$ when $k\geq 1$.
By Theorem \ref{ThLN16-BalBip}, $G$ is Hamiltonian or $G\subseteq B_n^k$. By Theorem 1.8,
$G$ is Hamiltonian or $G=B_n^k$.
The proof is complete. {\hfill$\Box$}

\section{Concluding remarks}
We suggest the following general problems.

\begin{problem}\label{Prob1}
Let $G$ be a connected graph of order $n$. Let
$s$ be an integer with $s\geq 1$.
Suppose that $\lambda(G)>\lambda(K_1\vee (K_s+K_{n-s-1}))$,
where $n$ is sufficiently large compared to $s$.
Does $G$ contain a $C_{n-s+1}$?
\end{problem}

\begin{problem}\label{Prob2}
Let $G$ be a connected graph of order $n$. Let
$s$ be an integer with $s\geq 1$.
Suppose that $q(G)>q(K_1\vee (K_s+K_{n-s-1}))$,
where $n$ is sufficiently large compared to $s$.
Does $G$ contain a $C_{n-s+1}$?
\end{problem}

Affirmative answers to these problems will give tight spectral
conditions for the existence of cycle $C_l$, where $l$ is large.
One can easily find that Theorems \ref{ThFN} and \ref{Th(n-1)cycleAQ}
give affirmative solutions to Problem \ref{Prob1} for the cases $s=1$
and $s=2$, respectively. Theorem \ref{Th(n-1)cycleAQ} solves
Problem \ref{Prob2} when $s=2$.

Moreover, we can also consider spectral conditions for consecutive cycles.
In this spirit, Theorems \ref{ThFN} and \ref{Th(n-1)cycleAQ} can be extended as follows,
respectively.

\begin{theorem}\label{ThAConsective}
Let $G$ be a graph of order $n\geq 5$. If $\lambda(G)>n-2$, then $G$ is pancyclic
unless $G=N^1_n$.
\end{theorem}

\begin{theorem}\label{ThQConsective}
Let $G$ be a graph of order $n\geq 15$.
If $\lambda(G)>\lambda(K_1\vee (K_2+K_{n-3}))$
or $q(G)> q(K_1\vee (K_2+K_{n-3}))$,
then $G$ contains a cycle $C_l$ for each
$l$ such that $3\leq l\leq n-1$.
\end{theorem}

The main ingredient of the proofs comes from a classical theorem proved
by Woodall \cite{W72}.

\begin{theorem}[\rm {Woodall \cite{W72}}]\label{ThBW}
Let $G$ be a graph of order $n\geq 2k+3$, where $k\geq 0$ is an integer. If
$$
e(G)\geq \binom{n-k-1}{2}+\binom{k+2}{2}+1,
$$
then $G$ contains a $C_l$ for
each $l$ such that $3\leq l\leq n-k$.
\end{theorem}

\noindent
{\bf Proof of Theorem \ref{ThAConsective}.}
Since $\Delta(G)\geq \lambda(G)>n-2$, we obtain $\Delta(G)\geq n-1$, which implies that $G$
is connected. By Theorem \ref{ThHong}, we get $\sqrt{2e(G)-n+1}>n-2$. We infer that
$2e(G)\geq n^2-3n+4\geq n^2-5n+14$ for $n\geq 5$. This implies that $e(G)\geq \binom{n-2}{2}+\binom{3}{2}+1$
for $n\geq 5$. By Theorem \ref{ThBW}, $G$ contains all cycles $C_l$,
where $3\leq l\leq n-1$. By Theorem \ref{ThFN}, $G$ contains a Hamilton cycle or $G=N^1_n$.
So $G$ contains all cycles of length from 3 to $n$. The proof is complete.  {\hfill$\Box$}

\noindent
{\bf Proof of Theorem \ref{ThQConsective}.}
If $\lambda(G)>\lambda(K_1\vee (K_2+K_{n-3}))$, then $\Delta(G)\geq \lambda(G)>n-3$.
If $q(G)> q(K_1\vee (K_2+K_{n-3}))$, then $2\Delta(G)\geq q(G)>2(n-3)$.
In each case, we deduce $\Delta(G)>(n-3)$,
which implies that $\Delta(G)\geq n-2$. It follows that $G$
contains at most one isolated vertex. If $G$ is connected, then
by Theorem \ref{ThHong}, we get $\sqrt{2e(G)-n+1}>n-3$. We infer that
$2e(G)\geq n^2-5n+8$. If $G$ is not connected, let $v$ be the isolated vertex and $G':=G-v$,
then we have $\sqrt{2e(G')-v(G')+1}>n-3$, that is $2e(G)\geq n^2-5n+8$. Since
$n^2-5n+8\geq n^2-6n+23$ for $n\geq 15$, this implies that $e(G)\geq \binom{n-3}{2}+\binom{4}{2}+1$
for $n\geq 15$.
By Theorem \ref{ThBW}, $G$ contains all cycles $C_l$,
where $3\leq l\leq n-2$. By Corollary \ref{Cor1}, $G$ contains a $C_{n-1}$.
The proof is complete.  {\hfill$\Box$}
\section*{Acknowledgements}
Jun Ge is supported by NSFC (Nos. 11701401 and 11626163). Bo Ning is
supported by NSFC (Nos. 11601379 and 11771141). The revised version was finished
when the authors were visiting Professor Fengming
Dong at NTU, Singapore. The authors are very grateful
to all his help, enthusiasm and encouragement during
their visits. The authors would
like to express their deep gratitude to the referee
for his/her suggestions that considerably
improve the quality of the paper.

\end{document}